\documentclass[10pt,letterpaper]{article}
\usepackage[utf8]{inputenc}
\usepackage[scanall]{psfrag}
\usepackage{graphicx}
\usepackage{epstopdf}
\usepackage{tikz}
\usetikzlibrary{shapes,arrows,chains}
\usepackage[utf8]{inputenc}
\usepackage{mathtools}
\usepackage{multirow}
\usepackage{color}
\usepackage{here}
\usepackage{array}
\usepackage{scrextend}
\usepackage{booktabs}
\usepackage{tabularx}
\usepackage{listings}
\usepackage{color}
\usepackage[ruled,vlined]{algorithm2e}
\usepackage{adjustbox}
\usepackage{setspace}
\usepackage[breaklinks]{hyperref}
\usepackage{fancyhdr}
\usepackage[noadjust]{cite}

\newlength{\bibitemsep}
\newlength{\bibparskip}
\let\oldthebibliography\thebibliography
\renewcommand\thebibliography[1]{%
  \oldthebibliography{#1}%
  \setlength{\parskip}{\bibitemsep}%
  \setlength{\itemsep}{\bibparskip}%
}

\usepackage[margin=1in]{geometry}
\usepackage{censor}

\allowdisplaybreaks

\newcolumntype{D}{>{\raggedleft\arraybackslash}X}
\newcolumntype{S}{>{\scriptsize\raggedleft\arraybackslash}X}

\usetikzlibrary{matrix,shapes,arrows,positioning,chains}
\usepackage{xspace}
\usepackage{float}

\graphicspath{ {./Figures/} }

\newcommand{\mV}{\vspace{1mm}}

\newcommand{\mmmV}{\vspace{5mm}}

\newcommand{\T}[1]{\text{#1}}
\newcommand{\mi}{\T{--}}

\newcolumntype{Y}{>{\small\centering\arraybackslash}X}

\fancypagestyle{mainpage}{\lhead{}\rhead{}\chead{Skolfield, Ramirez-Vergara, and Escobedo}}

\begin{document}
\fancypagestyle{firstpage}{
\lhead{Proceedings of the 2021 IISE Annual Conference \\ 
A. Ghate, K. Krishnaiyer, K. Paynabar, eds.}}
\renewcommand{\headrulewidth}{0pt}
\thispagestyle{firstpage}

\begin{center}
    \fontsize{16pt}{19.2pt}\selectfont{\textbf{Transmission and Capacity Expansion Planning Against Rising Temperatures: A Case Study in Arizona}}
    
    \mV 
    
    
    \fontsize{12pt}{14.4pt}\selectfont{\textbf{J. Kyle Skolfield and Adolfo R. Escobedo \\ Arizona State University \\ Tempe, Arizona}}
    
    \fontsize{12pt}{14.4pt}\selectfont{\textbf{Jose Ramirez-Vergara \\ Purdue University \\ West Layfayette, Indiana}}
    
    
    \fontsize{12pt}{14.4pt}\selectfont{\textbf{Abstract}}
\end{center}
The stable and efficient operation of the transmission network is fundamental to the power system’s ability to deliver electricity reliably and cheaply. As average temperatures continue to rise, the ability of the transmission network to meet demand is diminished. Higher temperatures lead to congestion by reducing thermal limits of lines while simultaneously reducing generation potential. Due to prohibitive costs and limited real estate for building new lines, it is necessary to consider capacity expansion as well to improve the functioning and efficiency of the grid. Optimal control, however, requires many discrete choices, rendering fully accurate models intractable. Furthermore, temperature changes will impact different regions and climate differently. As such, it is necessary to model both temperature changes and transmission flows with high spatial resolution. This work proposes a case study of the transmission grid centered in Arizona, using a DC optimal power flow mathematical formulation to plan for future transmission expansion and capacity expansion to efficiently meet demand. The effects of rising temperatures on transmission and generation are modeled at the regional level. Several classes of valid inequalities are employed to speed up the solution process. Multiple experiments considering different temperature and demand trends are considered which include each of the above technologies.

\fontsize{12pt}{14.4pt}\selectfont{\textbf{Keywords}}\normalsize

Energy Systems, Transmission Expansion Planning, Optimization, Mathematical Programming, Valid Inequalities


\fontsize{12pt}{14.4pt}\selectfont{\textbf{1. Introduction\\}}\normalsize 
Environmental stressors are among the main causes of power system disturbances worldwide \cite{panteli2015influence}. High temperatures, for instance, can limit the transfer capability of transmission lines by increasing energy losses and line sagging \cite{panteli2015influence}. Rising temperatures also disrupt demand patterns, as they have been proved to be positively correlated with extreme temperatures during summertime \cite{miller2008climate}. Preparing power systems for adverse weather is a challenging task since the frequency of these events is increasing, but their precise occurrence in time and space is not known ahead of time.  Increases in power system robustness and reliability are among the most common long-term methods to guarantee the operation of the network during extreme circumstances. Generation Expansion Planning (GEP) and Transmission Expansion Planning (TEP) have been proposed as possible pathways to adapt the system to new conditions \cite{li2014electric}. The TEP problem seeks to reinforce the transmission network and provide a stable supply even under worst-case scenarios. TEP has both economic and engineering reliability objectives; this makes the problem a complex case of multi-objective optimization \cite{garcia2016dynamic}. A common approach to solving the TEP problem is the use of the DC optimal power flow (DCOPF) approximation. In this formulation, both economic dispatch and optimal power flow are considered in modeling the problem. The DC formulation offers a good approximation of the AC power flow for planning purposes, especially since it is faster and easier to solve \cite{wood2013power}.  The impacts of climate change, and specifically rising ambient temperature, have been studied from the generation and demand perspective in generation and transmission expansion problems \cite{hejeejo2017probabilistic, sathaye2013estimating, mcfarland2015impacts}. TEP has also been solved using a decentralized approach in which the electricity network was divided into regions to account for differences in demand and supply sides \cite{hariyanto2009decentralized, de2008transmission}. On the other hand, the effects of rising temperatures on transmission lines themselves were analyzed and estimated to account for capacity reduction in \cite{bartos2016impacts}.

This study proposes a DCOPF formulation with discrete transmission decisions and regional temperature considerations to solve the transmission expansion problem. The main contribution is the inclusion of a capacity reduction factor for the transmission lines, and the division of the electricity network by considering different climate regions to serve as a bridge between those studies solving TEP and those studying the impact of climate change on transmission networks. The analysis also includes a case study of the Arizona transmission network, calculating the results for 16 different discrete scenarios to account for differences in ambient temperature expectations across regions determined from an analysis of historical data. 


\pagestyle{mainpage}
\fontsize{12pt}{14.4pt}\selectfont{\textbf{2. Methodology\\}}\normalsize
For the purpose of this paper, the methodology first characterizes temperature-driven planning scenarios. This involves estimating lower and upper bounds for ambient temperature, dividing Arizona into climate regions, and determining the implications of rising temperature on the capacity of transmission lines. These data are used as inputs for the featured TEP model, which considers two types of transmission investment options: expansion via new lines and expanding the capacity of existing lines.  Data collection was performed using ambient temperature historical registers provided by the National Oceanographic and Atmospheric Administration (NOAA) and its National Centers for Environmental Information \cite{ncei}. This consists of 70 years' worth of data in daily increments, reporting maximum, minimum, and average daily temperature for selected locations across the state of Arizona.

\begin{figure}[!htb]
    \minipage{0.32\textwidth}
    \centering
    \includegraphics[width=.8\linewidth]{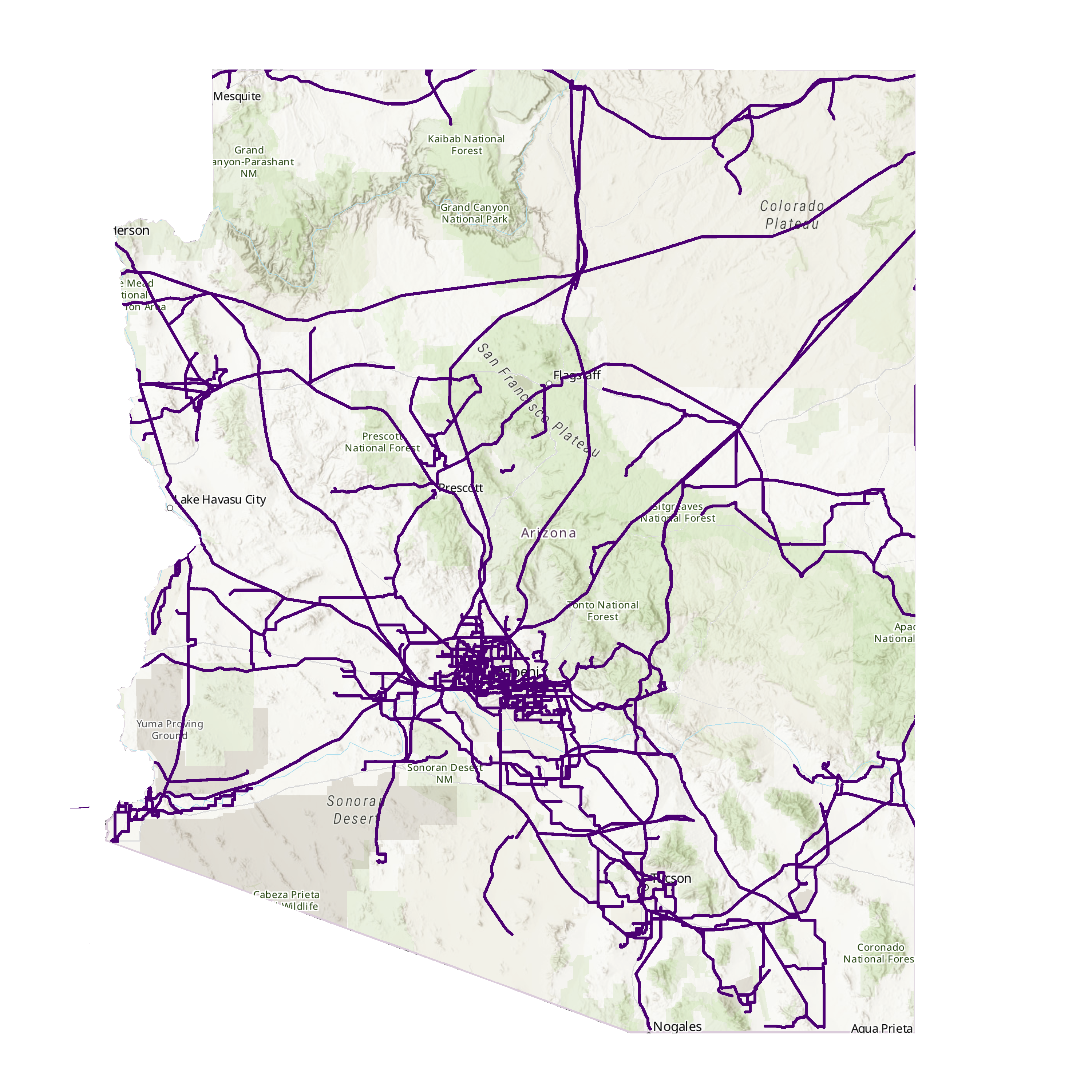}
    \caption{Transmission Network}
    \label{fig:trans_net}
    \endminipage\hfill
    \minipage{0.32\textwidth}
    \centering
    \includegraphics[width=.8\linewidth]{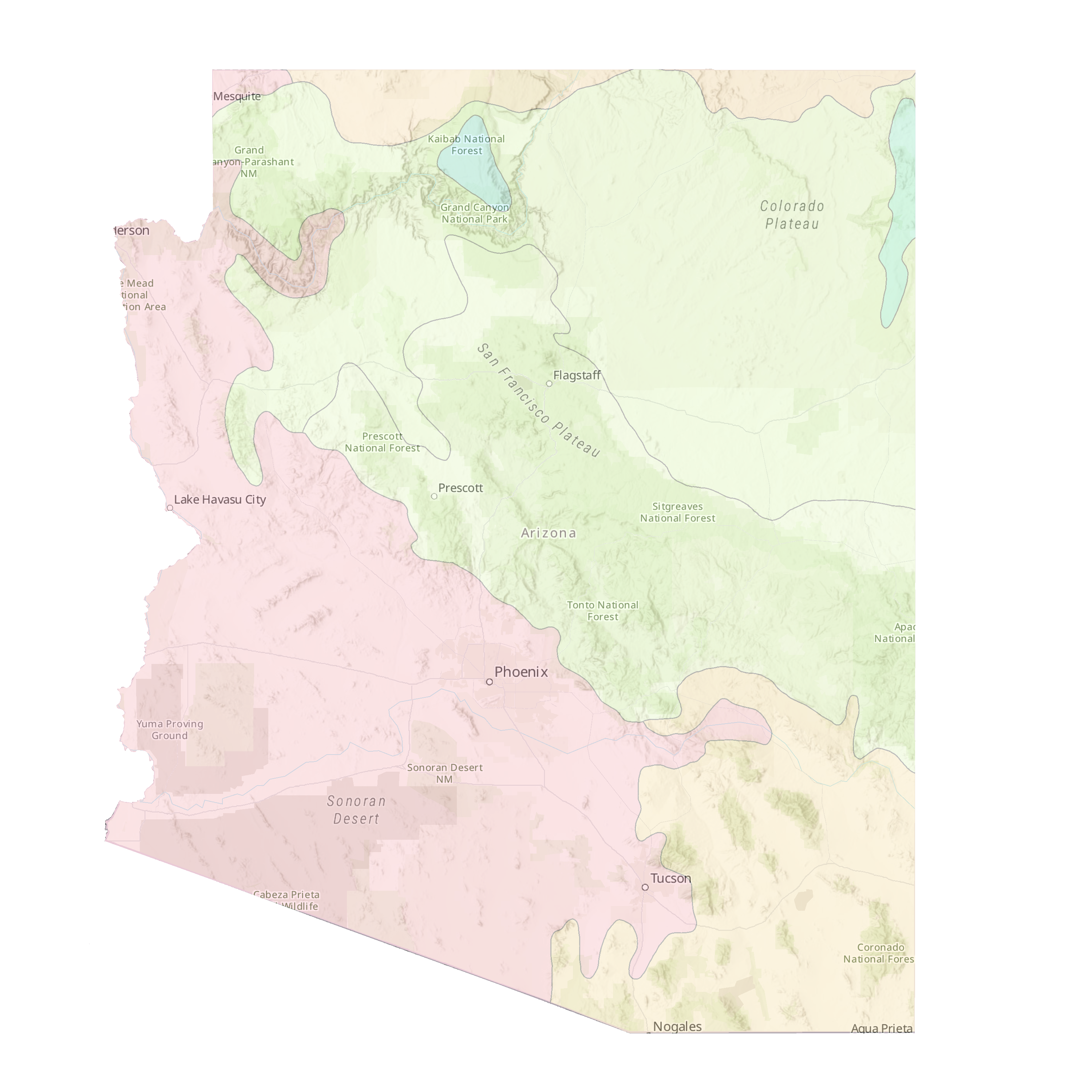}
    \caption{EPA Designated Ecoregions}
    \endminipage\hfill
    \minipage{0.32\textwidth}
    \centering
    \includegraphics[width=.8\linewidth]{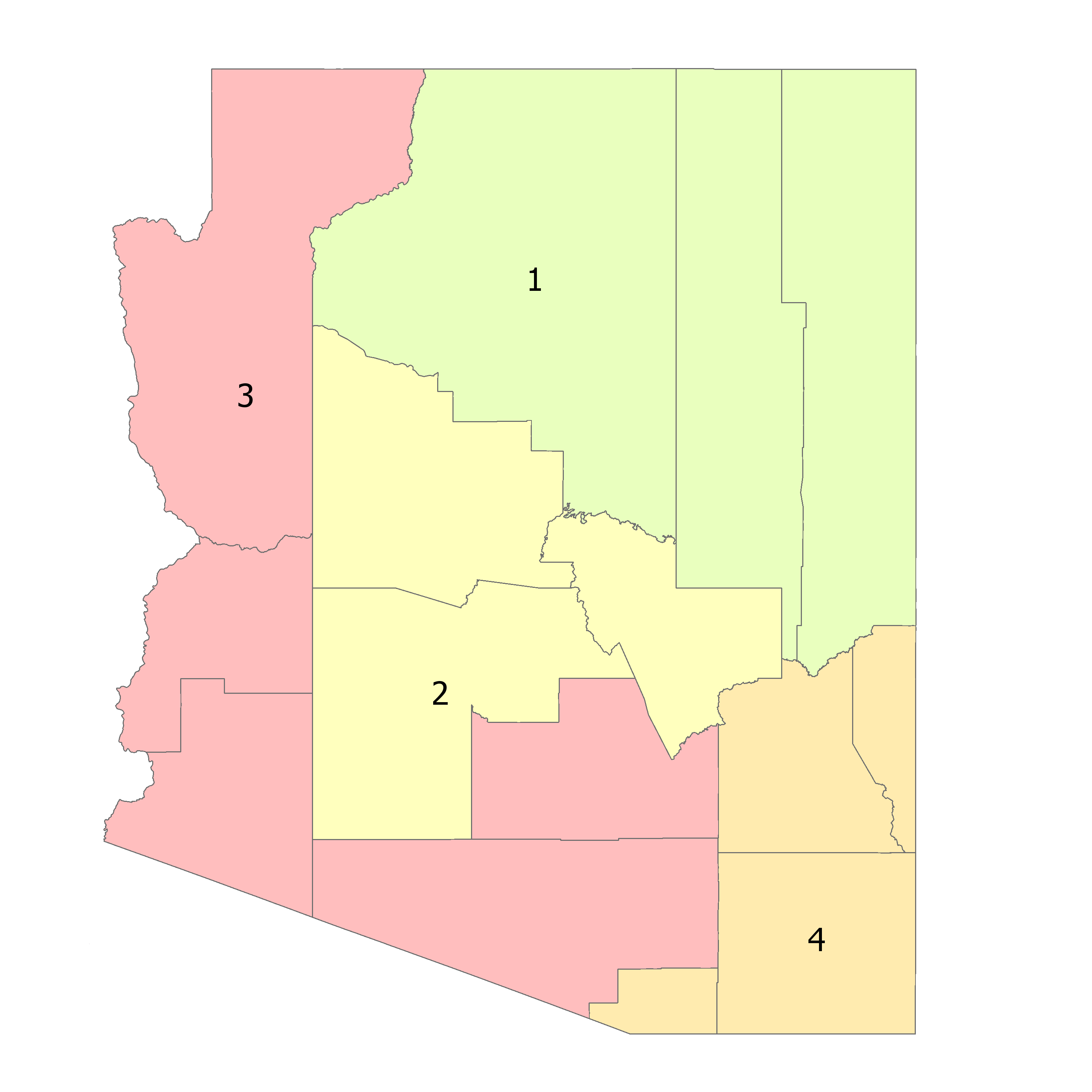}
    \caption{County Region Definitions}
    \label{fig:county}
    \endminipage\hfill
\end{figure}

\fontsize{10pt}{12pt}\selectfont{\textbf{2.1 Definition of Climate Regions\\}}\normalsize
The climate regions utilized in this work are based on the level II ecoregions defined by the Environmental Protection Agency; note the explanation for the development of these regions is provided in \cite{omernik2014ecoregions}.  Arizona was then divided into four major regions, with each county being assigned to a region based on the ecoregion which comprises the largest area within that county.  Figures \ref{fig:trans_net}-\ref{fig:county} demonstrate the transmission network, climate regions and corresponding county designations respectively.  

\mmmV

\fontsize{10pt}{12pt}\selectfont{\textbf{2.2 Estimation of Temperature Bounds\\}}\normalsize
Datasets incorporating daily temperature records for the period 1950-2019 were downloaded from NOAA for four representative urban centers (Phoenix, Tucson, Flagstaff, and Douglas). The 10 highest temperatures of each of these years were averaged to avoid outliers, and then plotted to analyze the behavior of maximum temperatures over time. The time series also provided enough data to perform regression analysis on the trend for peak annual temperature. This value corresponds to the expected maximum temperature for a given year. The highest maximum temperature per year was determined by selecting the data points above the mean trend line; such values and their associated years  used as inputs to another regression analysis which provided a linear equation for the overall maximum temperature trend. Both equations serve as a linear estimation for the future trend (30 years ahead), following the approach used in \cite{garcia2015multi}. These values were validated in comparison to the expected increase predicted by NOAA, often used as the reference forecast temperatures in other studies. NOAA forecasts temperature increases impacted by $CO_{2}$ emissions. The low emissions scenario projects increases in temperature of 2.7°F, 3.6°F, and 4.7°F by 2035, 2055 and 2085, respectively. The expected values due to high emissions are 3°F, 4.8°F, and 8°F for the same years \cite{Hayhoe2017}. By comparison, the projected trends computed from our regression analysis for Phoenix, AZ, suggest an increase of approximately  2.6°F by 2035, 3.60°F by 2055 and 5.1°F by 2085. These align closely with the intervals for low and high emissions scenarios provided by NOAA, hence validating their suitability for this study.

\fontsize{10pt}{12pt}\selectfont{\textbf{2.3 Impact of Temperature on Transmission Line Capacity\\}}\normalsize
The main consequence of rising temperatures on bulk power systems is the reduction in both transmission and generation capacity and the corresponding increase in demand \cite{bartos2016impacts, bartos2015impacts}. Considering the scope of this study, the estimation of the impact of rising temperatures is focused on electricity transmission. Following the Institute of Electrical and Electronics Engineers (IEEE) Standard 738-2006, the formula to calculate the reduction in the ampacity of the transmission lines also used in \cite{bartos2016impacts,abdalla2013weather,michiorri2015forecasting,jiang2016dispatching}, will be applied in this paper. The formula, derived from the energy balance equation, is as follows:
\begin{equation}
    I=\sqrt{\frac{\pi \cdot \overline{h} \cdot D \cdot (T_{cond} - T_{amb}) + \pi \cdot \epsilon \cdot \sigma \cdot D \cdot (T_{cond}^4 - T_{amb}^4) - \delta \cdot D \cdot a_{s}}{R \cdot T_{cond}}} \label{eqn:ampacity}
\end{equation}

Where $I$ is the fractional multiplier to the rated capacity of conductor (ampacity), $\overline{h}$ is the average heat transfer coefficient, $D$ is the diameter of the conductor, $T_{cond}$ and $T_{amb}$ are the average conductor temperature and ambient temperature, respectively. This first set of terms corresponds to losses due to convection. The second set of terms, corresponding to the loss due to radiation, includes $\epsilon$ which is the emissivity of conductor surface, $\sigma$ which is the Stefan-Boltzmann constant, and the diameter, plus the conductor temperature, and the ambient temperature. In the last set of terms, $\delta$ is the maximum solar radiation, $a_{s}$ is the absorptivity of the conductor surface and includes the diameter of the conductor, to account for the heat gain from solar radiation. The dividing term $R$ corresponds to the AC resistance of the conductor. 

\fontsize{12pt}{14.4pt}\selectfont{\textbf{3. Model Formulation\\}}\normalsize 
TEP with the DCOPF approximation can be formulated as a disjunctive, mixed-integer linear model.  The objective \eqref{eqn:htls_obj} is to minimize the joint cost of generation and investment in both new lines and capacity expansion on existing lines.  The full model is as follows.
\small 
\begin{flalign}
& \hspace{0pt}\displaystyle  \min  \sum_{(i,j) \in \Omega} \left( c_{ij}y_{ij} + h_{ij}z_{ij} \right) + \sum_{n \in B} \sigma c_n g_n & \label{eqn:htls_obj} \\
& s.t. \hspace{0pt}\displaystyle  \mathrlap{\sum_{(n,i)\in\Omega} \left(P_{ni}^{0} +  P_{ni} \right) -  \sum_{(i,n)\in\Omega} \left(P_{in}^{0} + P_{in} \right) + g_{n} = \gamma_r d_{n}} &\forall n \in r, r \in R \label{eqn:htls_balance} \\
& \hspace{0pt}\displaystyle - \eta_r \left( \overline{P}^0_{ij} - \overline{P}^1_{ij} z_{ij} \right) \leq P_{ij}^{0} \leq \eta_r \left( \overline{P}^0_{ij} + \overline{P}^1_{ij} z_{ij} \right) &\forall(i,j) \in \Omega \setminus \Psi, r \in R \label{eqn:htls_capacity}  \\
& \hspace{0pt}\displaystyle  - \eta_r \left(\overline{P}_{ij} y_{ij}\right) \leq P_{ij}  \leq \eta_r \left( \overline{P}_{ij} y_{ij}\right) &\forall(i,j) \in \Psi, r \in R \label{eqn:htls_tep_capacity} \\
& \hspace{0pt}\displaystyle  \frac{\mi 1}{b_{ij}}P_{ij}^{0} - (\theta_{i}-\theta_{j})= 0 &\forall(i,j) \in \Omega \setminus \Psi  \label{eqn:htls_flow} \\
& \hspace{0pt}\displaystyle \mathrlap{ \mi M_{ij}(1-y_{ij}) \leq \frac{\mi 1}{b_{ij}}P_{ij}-(\theta_{i}-\theta_{j}) \leq M_{ij}(1-y_{ij})} &\forall(i,j) \in \Psi \label{eqn:htls_tep_flow}  \\
& \hspace{0pt}\displaystyle  0 \leq g_{n} \leq \overline{g}_n   &\forall n \in B & \label{eqn:htls_gen} \\
& \hspace{0pt}\displaystyle  -\overline{\theta} \leq \theta_{i}-\theta_{j} \leq \overline{\theta} &\forall (i,j) \in \Omega  \label{eqn:htls_busangle} \\
& \hspace{0pt}\displaystyle  z_{ij} + y_{ij} \leq 1 &\forall (i,j) \in \Phi \label{eqn:htls_mutual}\\
& \hspace{0pt}\displaystyle  y_{ij}, z_{ij}  \in  \left\{  0,1  \right\}  &\forall(i,j) \in \Psi \forall(i,j) \in \Omega \setminus \Psi
\end{flalign}\normalsize
Throughout, $B$ is the set of all buses, $R$ is the set of regions, $\Omega$ is the set of all lines $(i,j)$, $\Psi$ is the set of candidate lines and $\Phi$ the set of expandable lines.  Constraint \eqref{eqn:htls_balance} enforces flow balance of $P_{ij}$, load $d_n$, and generation $g_n$, at each node $n$ in the transmission network.  The demand is scaled by the factor $\gamma_r$, which is set at the upper and lower limits (for high and low temperature increases, respectively) of predicted demand increase in \cite{bartos2016impacts}.  Constraints \eqref{eqn:htls_capacity} and \eqref{eqn:htls_tep_capacity} represent capacity limits $\overline{P}_{ij}$ on existing lines and candidate lines respectively.  The coefficient $\eta_r$ is the scaling factor induced by \eqref{eqn:ampacity} to account for per-region temperature change; specifically, $\eta_r$ is the value of $I$ obtained by substituting in the ambient temperature predicted for each region.  Constraints \eqref{eqn:htls_flow} and \eqref{eqn:htls_tep_flow} relate adjacent bus angles to power flow on the transmission line connecting them.  Constraint \eqref{eqn:htls_mutual} states that only existing lines can have their capacity expanded, whereas candidate lines can be built but not then expanded.  The remaining constraints represent domain limits on decision variables.  

\fontsize{10pt}{12pt}\selectfont{\textbf{3.1 Valid Inequalities\\}}\normalsize
TEP is known to be NP-Hard \cite{latorre}.  As such, a set of valid inequalities adapted from \cite{Skolfield2019} are incorporated into the solution algorithm to reduce computational time.  The model could not be solved to optimality within 24 hours without the valid inequalities; it was solved for all instances under this time limit after their inclusion.  The full statement of the valid inequalities is given by 
\small
\begin{align*}
    \vert \theta_n - \theta_m \vert \leq CR^r\left(\rho \setminus \left(i,j\right)\right) + \left( \overline{CR^r(\rho)} - CR^r\left(\rho \setminus \left(i,j\right)\right) \right) \left( N_e(\rho) - \sum_{(l,k) \in \rho \cap \Psi} y_{lk} \right) +
    \left( {CR^r(\rho)} - CR^r\left(\rho \setminus \left(i,j\right)\right) \right) z_{ij}.
\end{align*}\normalsize
Here, $CR^r$ is the sum of the capacity-reactance product of each line in the path $\rho$, using the expanded capacity for all lines that can be reconductored.  $N_e(\rho)$ counts the number of candidate lines on $\rho$.

\fontsize{12pt}{14.4pt}\selectfont{\textbf{4. Instance Generation\\}}\normalsize
The case study presented in this paper is an approximation of the transmission grid contained within the borders of Arizona; the basic description of this network is provided by the U.S. Energy Information Administration.  Transmission lines rated at 69kV and above are included in the network as are all power plants rated to produce two or more MW/Hr.  Additional corridors for transmission expansion are approximated by connecting high generation areas to high demand areas and connecting substations with few adjacent lines to more dense areas of the grid.  All listed substations are used as buses to connect transmission lines and meet aggregated demand for nearby areas.  Hourly demand for the test instance is based on peak summer demand for the full state.  As the state level is the smallest resolution data available for load, this value is disaggregated by assigning load to substations according to the relative population in the nearest census block.  Generation costs are approximated by the total statewide costs of generating each category of plant (natural gas, coal, petroleum, hydro, wind, and solar), divided according to the rated MW of the corresponding plant.  Transmission expansion and capacity expansion costs are approximated based on the rated voltage and length of lines, using cost estimates from \cite{MISO}.

\fontsize{12pt}{14.4pt}\selectfont{\textbf{5. Results\\}}\normalsize
A MILP approximation of the AZ transmission network is solved for 16 scenarios: each climate zone of the network is projected to have either a large or small increase (e.g., $2\mi5^{\circ}$F for Tucson) in summer peak temperature independent from each other zone.  These scenarios are encoded by an ordered quadruplet in which the temperature increase of region $i$ is indicated by an H in position $i$ if it is a large increase, and an L in the same position if a small increase is projected.  This model is solved in Gurobi 9.0.2, and the non-zero expansion variables (and their associated objective costs) are tallied for each scenario.  The results of these experiments are summarized in Table \ref{tab:my-table}, with costs annualized.

\begin{table}[!hbtp]
\caption{Per Annum Projected Cost}
\label{tab:my-table}
\begin{tabularx}{\textwidth}{|Y|Y|Y|Y|Y|Y|Y|Y|} \hline 
\textbf{Scenario} & \textbf{New Lines Built} & \textbf{Cap. Exp. Built} & \textbf{New Line Cost} & \textbf{Cap. Exp. Cost} & \textbf{Total Exp. Cost} & \textbf{\;\;\;Gen.\newline Cost} & \textbf{\;\;Total\newline Cost} \\\hline\hline
 L,L,L,L           & 80                & 17               & \$  14.40B    & \$  2.79B    & \$     17.19B        & \$    7.07B      & \$  24.26B  \\
 L,L,L,H           & 83                & 19               & \$  14.94B    & \$  3.15B    & \$     18.09B        & \$    6.39B      & \$  24.48B  \\
 L,L,H,L           & 81                & 23               & \$  14.58B    & \$  3.96B    & \$     18.54B        & \$  11.27B       & \$  29.81B  \\
 L,H,L,L           & 89                & 17               & \$  16.02B    & \$  2.85B    & \$     18.87B        & \$    6.06B      & \$  24.93B  \\
 H,L,L,L           & 76                & 18               & \$  13.68B    & \$  2.94B    & \$     16.62B        & \$    7.47B      & \$  24.09B  \\
 L,L,H,H           & 83                & 25               & \$  14.94B    & \$  4.24B    & \$     19.18B        & \$  10.97B       & \$  30.15B  \\
 L,H,L,H           & 84                & 18               & \$  15.12B    & \$  3.00B    & \$     18.12B        & \$    6.77B      & \$  24.89B  \\
 H,L,L,H           & 85                & 20               & \$  15.30B    & \$  3.24B    & \$     18.54B        & \$    6.00B      & \$  24.54B  \\
 L,H,H,L           & 70                & 24               & \$  12.60B    & \$  4.09B    & \$     16.69B        & \$  13.21B       & \$  29.90B  \\
 H,L,H,L           & 88                & 24               & \$  15.84B    & \$  4.09B    & \$     19.93B        & \$  10.22B       & \$  30.15B  \\
 H,H,L,L           & 81                & 21               & \$  14.58B    & \$  3.41B    & \$     17.99B        & \$    6.98B      & \$  24.97B  \\
 L,H,H,H           & 83                & 22               & \$  14.94B    & \$  3.73B    & \$     18.67B        & \$  11.46B       & \$  30.13B  \\
 H,L,H,H           & 85                & 22               & \$  15.30B    & \$  3.79B    & \$     19.09B        & \$  10.70B       & \$  29.79B  \\
 H,H,L,H           & 83                & 18               & \$  14.94B    & \$  3.00B    & \$     17.94B        & \$    7.32B      & \$  25.26B  \\
 H,H,H,L           & 89                & 22               & \$  16.02B    & \$  3.73B    & \$     19.75B        & \$  10.16B       & \$  29.91B  \\
 H,H,H,H           & 78                & 22               & \$  14.04B    & \$  3.73B    & \$     17.77B        & \$  12.50B       & \$  30.24B \\ \hline 
\end{tabularx}
\end{table}

Total costs increase by over 25\% when comparing the scenario in which all regions experience small temperature gains to the scenario in which all regions experience large temperature gains.  Since the difference in temperature changes for all regions is less than 3\%, this represents a superlinear relationship between temperature increases and associated power system costs.  Furthermore, larger increases in temperatures have different effects depending on the regions in which they occur.  For example, when regions 1 and 4 experience large temperature increases, the associated generation costs are much lower than when either region 2 or 3 experiences larger temperature increases.  Since regions 2 and 3 contain Phoenix and Tucson respectively, the most populous metropolitan areas in the state by large margins, this is consistent with expectations.  Generation costs associated with higher temperatures in these regions are nearly twice as much as those associated with higher temperatures in the less populous regions 1 and 4.  

Comparing the relative volume and costs associated with new lines rather than capacity expansion reveals further features of this network.  When region 3 (containing Tucson) experiences larger temperature increases, more lines are reinforced with higher capacity than in any other region.  This suggest that Tucson has a robust infrastructure in place when considering transmission lines and generation, so the relatively less expensive option of expanding the capacity of existing lines can accommodate a large amount of increased demand in the region.  In contrast, when region 2 (containing Phoenix) experiences larger temperature increases, new lines are built at a higher rate (and correspondingly, cost) than for any other region.  The increase in demand associated from higher temperatures in this region does not cause generation costs to increase very much, but it does require a large number of new lines to be built to meet this demand.  This is consistent with a less robust transmission infrastructure in this region compared to region 3, but better access to current and large sources of generation.

\fontsize{12pt}{14.4pt}\selectfont{\textbf{6. Conclusions\\}}\normalsize
Significant temperature increases are expected by mid century: the question is the magnitude and distribution of these increases.  This study considers several scenarios of temperature increase, distributed across both magnitude and location, and their effects on state-level transmission networks.  A test case is built from U.S. Government provided data and publicly available data on the existing generation and transmission assets within Arizona; temperature scenarios are similarly designed based on historical data and designated ecoregions.  The optimal expected annual cost of operation the power network, as well as costs associated solely with generation, transmission expansion, and capacity expansion is projected.  Based on the distribution of these costs, it can be seen that regardless of the nature of temperature increases, a significant investment in transmission infrastructure is required by mid century.  Furthermore, conditions associated with the existing network will cause the total cost to depend -- to varying degrees -- on transmission expansion, capacity expansion, or generation.  Large increases in temperatures over urban areas correspond to even larger increases in generation costs in those regions.  The largest overall cost differences are associated with the degree of temperature increase in these urban regions.  There are also network features which are distinct to each such region that dictate whether future demand can be met primarily with less expensive capacity expansion via reconductoring or require significantly more costly transmission expansion to build new overhead lines.  These results suggest that further analysis is worth performing, including cost changes due to generation or substation investments.  The demand projections are also limited in scope to a fixed percentage increase based on the regional temperature changes.  However, the current work demonstrates that such analysis is valuable.  Further work is necessary, including an analysis of some subnetworks of the Arizona power system with more varied scenarios, in order to fully understand what network features suggest investment in certain classes of transmission asset.  The joint optimization of both generation and transmission, especially with significant state- and nationally-mandated plans to expand renewable generation, also remains a question of interest.  

\bibliographystyle{IEEEtran}
\bibliography{references}

\end{document}